\documentclass[11pt] {article}      %
\usepackage{amsmath,amssymb, amsthm, eucal}  
\usepackage{graphicx, amscd}

\usepackage{txfonts}                  
\usepackage[T1]{fontenc}     

\usepackage{float}                                       
\usepackage[hang,small,sc]{caption}             

\theoremstyle{plain}
\newtheorem{theorem}{Theorem}[section]            
\newtheorem{proposition}[theorem]{Proposition}  
\theoremstyle{definition}
\newtheorem{definition}[theorem]{Definition}    
\newtheorem{remark}[theorem]{Remark}	
\newtheorem*{xremark}{Remark}		

\numberwithin{theorem}{section}
\numberwithin{equation}{section}
\numberwithin{figure}{section}

\usepackage{sidecap}    

\usepackage{titlesec}                                               
\titlelabel{\thetitle.\ }                                                
\titleformat*{\section}{\fontsize{14pt}{14pt} \bf}        
\usepackage{fancyhdr}
\begin{document}

\title{A porism concerning cyclic quadrilaterals}

\author{Jerzy Kocik 
    \\ \small Department of Mathematics, Southern Illinois University, Carbondale, IL62901
   \\ \small jkocik{@}siu.edu
}

\date{}

\maketitle

\begin{abstract}
\noindent
We present a geometric theorem on a porism about cyclic quadrilaterals, namely the existence of an
infinite number of cyclic quadrilaterals through four fixed collinear points once one exists. 
Also, a technique of
proving such properties with the use of pseudounitary traceless matrices is presented. 
A similar property holds for general quadrics as well as the circle.
\footnote{Published in {\it Geometry}, Volume 2013 (Jun 2013), Article ID 483727.}
\end{abstract}

\section{Introduction} 

\begin{figure}[H]
\centering
\includegraphics[scale=.9]{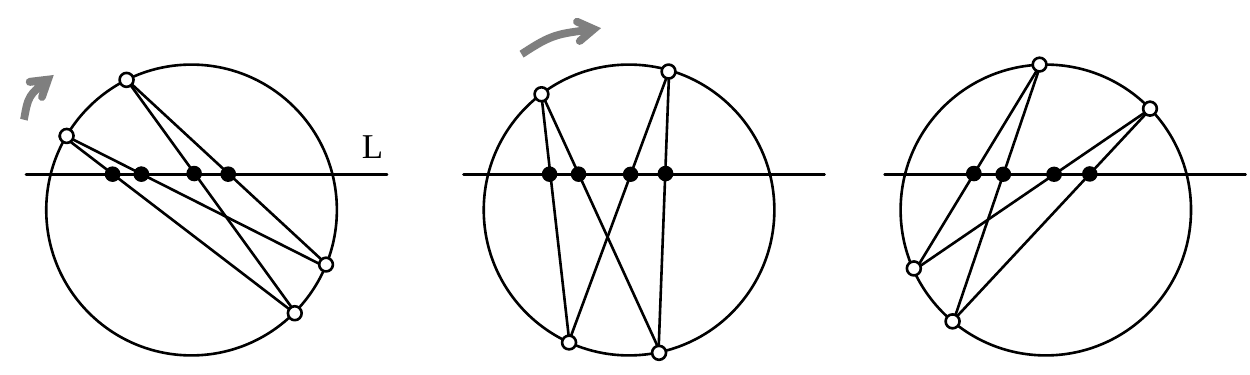}
\caption{\small  Porism on cocyclic quadrilaterals through collinear points}
\label{fig:porism}
\end{figure}

Inscribe a butterfly-like quadrilateral in a circle and draw a line $L$, see Figure \ref{fig:porism}. 
The sides of the quadrilateral will cut the line at four (not necessarily distinct) points.
 It turns out that as we  continuously deform the inscribed quadrilateral, the points of intersection remain invariant. 

\begin{figure}[h]
\centering
\includegraphics[scale=.8]{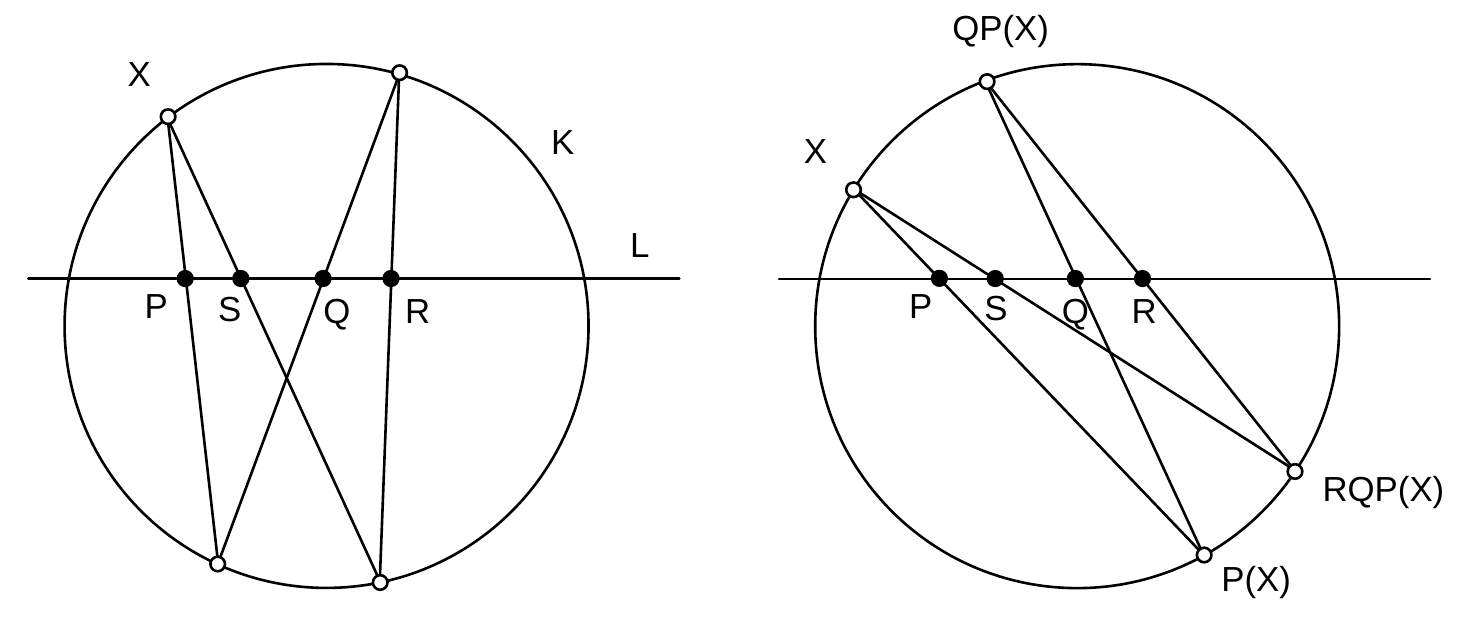}
\caption{\small Porism: details}
\label{fig:detail}
\end{figure}

More precisely, think of the quadrilateral as a path from vertex $X$ through collinear points $P$, $R$, $S$
and $Q$ back to $X$ (see Figure \ref{fig:detail}, left).  If we redraw the path starting from another point on the
circle but passing through the same points on the line in the same order, the path closes to form an inscribed polygon, 
that is, we shall arrive at the starting point. 

\begin{SCfigure}[2.1][h]   
\centering
\includegraphics[scale=.8]{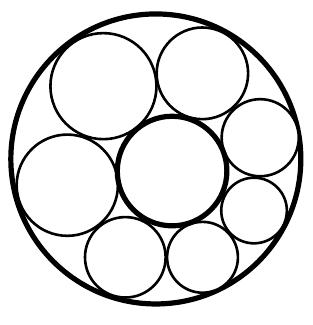}
\caption{\small Steiner porism. The chain may be redrawn starting from arbitrary position of the starting circle.}
\label{fig:steiner}
\end{SCfigure}

In spirit, this startling property is similar to Steiner's famous porism  \cite{CG, W1}, which states that once we find two
circles, one inner to the other, such that a closed chain of neighbor-wise tangent circles inscribed in the region
between them is possible, then an infinite number of such inscribed chains exist
(Figure \ref{fig:steiner}). 
One may set the initial circle at any position and the chain will close with tangency. 
Yet another geometric phenomenon in the same category is Poncelet's porism \cite{BB, Po, W2}. 


The property for the cyclic quadrilateral described at the outset may be restated similarly: if four points on a line 
admit a cyclic quadrilateral then an infinite number of such quadrilaterals inscribed in the same circle exist. 
Hence the term {\it porism} is justified.  
\\

In the next sections we restate the theorem and define a map of reversion through point, which may be represented 
by pseudo-unitary matrices.  The technique developed allows one to prove the theorem as well as 
a diagrammatic representation of the relativistic addition of velocities, presented elsewhere \cite{jkd}.  
A slight modification to arbitrary two-dimensional Clifford algebras allows one to modify the result to hold 
for hyperbolas and provides a geometric realization of trigonometric tangent-like addition. 
\\

\section{The main result}

\begin{SCfigure}[2.2][h]
\centering
\includegraphics[width=0.31\textwidth]{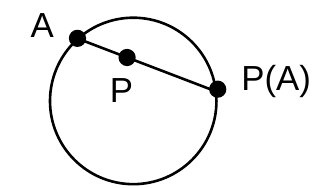}
\caption{\small Reversion of $A$ through $P$ }
\label{fig:reversion}
\end{SCfigure}

Let us present the result more formally.  Reversion of a point $A$ on a circle through point $P$  gives
point $\mathbf P (A)$ on the circle such that points $A$, $P$ and $\mathbf P(A)$ are
collinear (see Figure \ref{fig:reversion}). \ More precisely:

\begin{definition}
Given a circle $K$, a {\bf reversion} through point $P\not\in K$ is a map $\mathbf P$:
$K \rightarrow K,  \ A \mapsto \mathbf P(A)$, such that points
$A, P, \mathbf P(A)$ are collinear, and $\mathbf P(A) \neq A$. 
If $P \in K$, then for any $A\in K$  we define $\mathbf P(A) = P$. 
\end{definition}
%
%

By bold $\mathbf P$ we denote the reversion map defined by point $P$ (not bold).
Clearly, reversion through $P \not\in K$ is an involution,  $\mathbf P^2= \hbox{id}$, and therefore invertible.
The product of reversions is in general not commutative, 
$\mathbf{PQ} \neq \mathbf Q\mathbf P$.  
The main result may now be given a concise form:

\begin{theorem}
\label{thm:main}
Let $P, Q,R,S$ be collinear points, $K$ a circle, and
$\mathbf Z = \mathbf {SRQP}$ the composition of the corresponding reversions.
Then
$$
  \hbox{if} \qquad \exists X\in K \quad \mathbf Z(X) = X
                \qquad\hbox{then}\qquad
   \forall X\in K \quad \mathbf Z(X) = X
$$
\end{theorem}

The quadrilateral $[X, \mathbf{P} (X), \mathbf{QP}(X),\mathbf{RQP}(X)]$ may be
viewed as a member of a family parametrized by $X$. 
As $X$ moves along the circle, the quadrilateral continuously changes its
form while the points of intersection with the line remain invariant. 
The quadrilateral appears to be ``rotating'' with
points $P, Q, R, S$ playing the role of ``axes'' of this motion.
For an animated demonstration, see \cite{jkgeo}.

Here is an equivalent version of Theorem \ref{thm:main}:

\begin{theorem}
\label{thm:equiv}
Let $K$ be a circle and $L$ a line with three points $P$, $Q$, $R$.
Let $X$ be a point on the circle. The point $S$, an intersection of line $L$ with
line $[\mathbf{RQP}(X),\; X]$, does not depend on $X\in K$.
\end{theorem}

In other words, for any three collinear points $P,Q,R\in L$, there exists a
unique point $S$ on $L$ such that \textit{any} cyclic quadrilateral inscribed in $K$ passing through
$R,P,Q$ (in the same order) must also pass through $S$. 
\\

A more general statement holds:

\begin{theorem}
\label{thm:iff}
The composition of three point reversions $\mathbf{PQR}$ of a circle is a reversion if and only if
points $\mathbf P$, $\mathbf Q$ and  $\mathbf R$ are collinear. 
\end{theorem}


\begin{xremark}  
The figures present quadrilaterals in the butterfly shape for convenience.  
Clearly, they can have untwisted shape and also
the points of intersection can lie outside the circle and 
the line does not need to intersect the circle.
\end{xremark}

\newpage
\section{Reversion calculus -- matrix representation}

In order to prove the theorem we develop a technique that uses complex numbers and matrices. 
Interpret each point $P$ as a complex number $p\in \mathbb C$.
Without loss of generality, we shall assume that $K$ is the unit circle,  $K = \{ z\in \mathbb C\, : \, |z|^2 = 1\}$. 
Complex conjugation is denoted in two ways. Here is our main tool:

\begin{theorem}
\label{thm:moebius}
Reversion in the unit circle $K$ through point $p\in \mathbb C$ corresponds to a M\"obius transformation:
\begin{equation}
\label{eq:moebius}
\hat P: \ z\quad  \longrightarrow \quad  z' = \left[\begin{array}{rr}          1 & - p \\ 
                                                                                          \bar p & -1  
                                                             \end{array}\right] 
         \cdot z  = \frac{z-p}{\bar pz-1}   
\end{equation}
\end{theorem}

\noindent {\bf Proof:} 
 First, we check that reversions leave the unit circle invariant,  $|z|^2 = 1 \Rightarrow  |z'|^2 =1$ :
$$
|z'|^2  = \frac{z-p}{\bar p z - 1} \left( \frac{z-p}{\bar p z - 1}\right)^\ast
            = \frac{|z|^2-z\bar p - \bar p z +|p|^2}{|p|^2|z|^2 - \bar pz -\bar zp +1}
            = \frac{1 -  z\bar p - \bar p z +|p|^2}{|p|^2 - \bar pz -\bar zp +1} = 1
$$
Next, we show that for any $z\in K$, points $(z, p, z' )$ are collinear. 
We may establish that by checking that $(z' - p)$ differs from $(z - p)$ only by scaling by a real number. 
Indeed, take their ratio:
$$
\begin{aligned}
\frac{z' - p}{z-p} &=   \frac{    \frac{z-p}{\bar p z - 1} - p  }{ z-p}
                                       =   \frac{(z-p) - p(\bar pz - 1) }{(z-p)(\bar pz-1)}  \\
                              &= \frac{1-|p|^2}{\bar p z + p\bar z - |p|^2 - 1}
                                       = \frac{1-|p|^2}{2\hbox{Re}(\bar p z) - |p|^2 - 1}
                                        \ \in \mathbb R
\end{aligned}
$$ 
where, to get the first expression in the second line, we multiplied the numerator and the denominator 
by $\bar z$ and used the fact that $|z|^2=1$. 
It is easy (and not necessary) to check that  $\hat P ^2 = \hbox{id}$ (up to M\"obius equivalence).%
\hfill $\square$
\\

A more abstract definition of the matrix representation emerges.  Namely, denote the pseudo-unitary group of two-by-two complex matrices
$$
U(1,1) = \{ A\in \mathbb C\otimes\mathbb C \;|\;  A^\ast JA = J \}
$$
where  $J$ denotes the diagonal matrix $J = \hbox{diag} (1, -1)$, and the star ``$\ast$'' denotes usual Hermitian transpose. 
Let $\sim$ be the equivalence relation among matrices:  $A\sim B$ if there exists $\lambda\in \mathbb C$, 
$\lambda\neq 0$, such that $A=\lambda B$.  Then we have a projective version of the pseudo-unitary group:
$$
PU(1,1) = U(1,1)/\sim
$$
which is equivalent to its use for M\"obius transformations.  Now, any element of the group 
may be represented by a matrix up to a scale factor. The essence of Theorem~\ref{thm:main} is that 
reversions correspond to such matrices with vanishing trace, $\hbox{Tr}\, A = 0$. 
\\
\\
Note that we do not follow the tradition of normalizing the determinants, as is typically done 
for the representation of the modular group $PS\!L(2,\mathbb Z)$. 
\\

\section{Algebraic proof of the porism}

We shall now prove the main results. 
\\
\\
{\bf Proof of the main theorem:} 
Consider the  product of three consecutive reversions through points $P$, $Q$ and $R$, as represented by matrices
$$
M   =       
\left[\begin{array}{rr}          1 & - r \\ 
                                       \bar r & -1  
                                                             \end{array}\right]    
 \left[\begin{array}{rr}          1 & - q \\ 
                                        \bar q & -1  
                                                             \end{array}\right]    
 \left[\begin{array}{rr}          1 & - p \\ 
                                        \bar p & -1  
                                                             \end{array}\right]    
  = 
 \left[\begin{array}{rr}          1 - \bar pq +\bar qr -r\bar p   & - p +q - r + p\bar q r \\ 
                                                \bar p - \bar q +\bar r - \bar p q \bar r  &  -1 + p\bar q -q\bar r + \bar r p  
                                                             \end{array}\right]    
$$
Note that   $M_{12} = -\overline M_{21}$, but we still need to see whether $\hbox{Tr}\, M = 0$.   
Recall the assumption that $p$, $q$ and $r$ are collinear.  If $r = mp+nq$, where $m+n=1$, $m,n\in \mathbb R$,
we observe that the diagonal elements are real:
$$
 			M_{11}   =   1 -  m(\bar pq +p\bar q) + m|p|^2 -  n|q|^2  \in \mathbb R
$$
$$
 			M_{22}   =   - M_{11} .
$$
Dividing every entry by $M_{11}$ (M\"obius equivalence) we arrive at
$$
M =  
 \left[\begin{array}{cc}          1    & \frac{- p +q - r + p\bar q r}{1 - \bar pq -\bar qr + r\bar p} \\ 
                                          \frac{\bar p - \bar q +\bar r - \bar p q \bar r}{1 - \bar pq -\bar qr + r\bar p}  &  -1  
                                                             \end{array}\right]    
$$
which is visibly a matrix of reversion, defining uniquely the point of reversion:
$$
s  =    \frac{ p -q + r - p\bar q r}{1 - \bar pq -\bar qr + r\bar p}
$$
What remains is to verify that $s$ is collinear with $p$, and $q$ (and therefore $r$).  
This may be done algebraically by checking that $(s-p) / (r-q) \in\mathbb R$, but it also follows neatly by geometry:  
there are two points on the circle collinear with $P$, $Q$, and $R$, say $A$ and $B$.  
Since $\mathbf{RQP}(A) = B$ and therefore $\mathbf S(A)$ must be $B$, thus $S$ lies on the same line.  
We have proved that given a circle $K$, for any three points $P$, $Q$, $R$, on a line $L$ there exists 
a point $S$ on $L$ such that $\mathbf{SPQR} = \hbox{id}$,  or $\mathbf S=\mathbf{PQR}$ 
(because $\mathbf S^2= \hbox{id}$). The main theorem follows.  $\square$
\\

\begin{figure}[h]
\centering
\includegraphics[scale=.8]{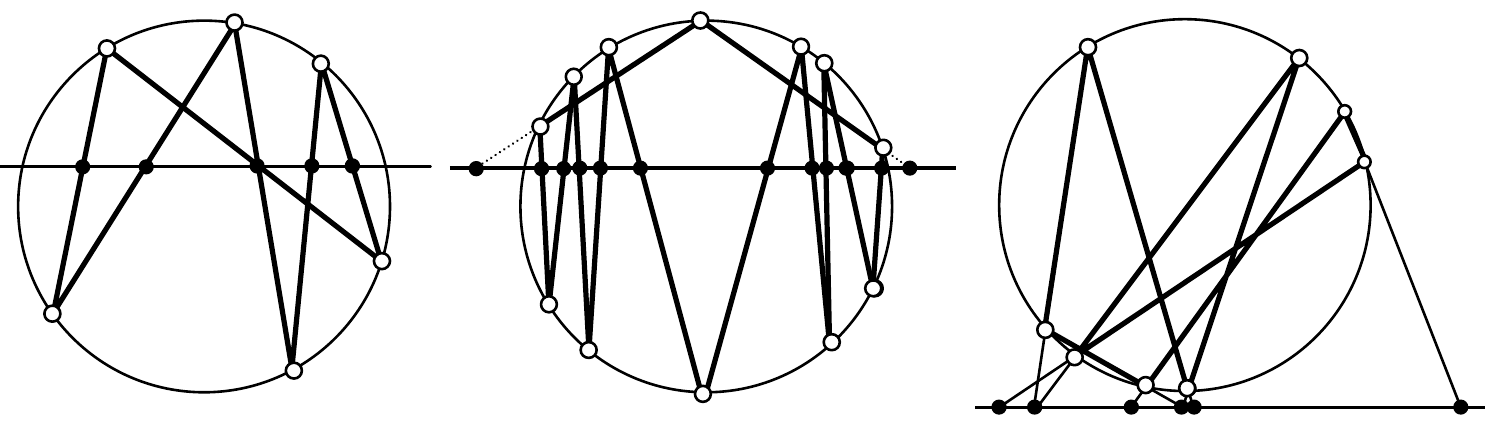}
\caption{\small A porism of cyclic polygons}
\label{fig:poly}
\end{figure}

The theorem generalizes to cyclic $n$-gons (see Figure \ref{fig:poly}):

\begin{theorem}
\label{thm:poly}  
Let  $\mathbf P = \mathbf P_{2n} \mathbf P_{2n-1}... \mathbf P_2 \mathbf P_1$  be the composition of 
point reversions of a circle $K$ for some collinear set of an even number of (not necessarily distinct) 
points $\mathbf P_1, ...,\mathbf P_{2n}$.  Then
$$
(\exists X\in K, \ \mathbf P(X) = X)  \quad \Rightarrow \quad (\forall X\in K , \ \mathbf P(X)=X) \; .
$$
\end{theorem}

\noindent
{\bf Proof:}  By Theorem \ref{thm:iff}, any three consecutive maps in the string $\mathbf P$  may be replaced by one.  
This reduces the string of an even number of reversions to a string of four, which is the case proven 
as Theorem \ref{thm:main}. 
Further reduction to two reversions would necessarily produce $\mathbf P_1\mathbf P_1= \hbox{id}$. 
\hfill $\square$

\begin{remark}
Theorem \ref{thm:poly} suggests a concept of a ``ternary semigroup of points'' on a (compactified) line, 
where the product requires three elements to produce one: 
$$
a,\;b,\;c \ \to\  abc
$$
The following may be considered as ``axioms'': 
\\

1. $(abc)de =  a(bcd)e = ab(cde)$	(associativity)

2. $abb = bba = a$   			(absorption of squares)

3. $abc = cba$    			(mirror inversion of triples)
\\
\\
The point reversions define a realization of this algebra.
\end{remark}

Yet another concept coming from these notes: Two pairs of points on a line are conjugate with respect to a circle if they belong to inscribed angles based on the same chord.  In other words: if  $\mathbf P\mathbf Q = \mathbf S\mathbf R$, 
or, if there exists an inscribed quadrilateral which intersects $L$ in $P$, $Q$, $R$, $S$. In matrix terms:
$$
M(q)M(p)  =  
 \left[\begin{array}{cc}          1  - q\bar p   &  q-p \\ 
                                               \bar q - \bar p   &  1 - \bar q p   
                                                             \end{array}\right]    
 =  \left[\begin{array}{cc}          1    &  \frac{q-p}{1  - q\bar p} \\ 
                                                \frac{\bar q - \bar p}{1  - q\bar p}   & \frac{ 1 - \bar q p}{1  - q\bar p}   
                                                             \end{array}\right]    
 =  \left[\begin{array}{cc}          1     &  \frac{q-p}{1  - q\bar p} \\ 
                                                \left(\frac{ q -  p}{1  - \bar q p}\right)^\ast   & \frac{ 1 - \bar q p}{1  - q\bar p}   
                                                             \end{array}\right]    \, ,
$$

$$
M(r)M(s)  =  
 \left[\begin{array}{cc}          1  - r\bar s   &  r-s \\ 
                                               \bar r - \bar s   &  1 - \bar r s   
                                                             \end{array}\right]    
 =  \left[\begin{array}{cc}          1    &  \frac{r-s}{1  - r\bar s} \\ 
                                                \frac{\bar r - \bar s}{1  - r\bar s}   & \frac{ 1 - \bar r s}{1  - r\bar s}   
                                                             \end{array}\right]    
 =  \left[\begin{array}{cc}          1     &  \frac{r-s}{1  - r\bar s} \\ 
                                                \left(\frac{ r -  s}{1  - \bar r s}\right)^\ast   & \frac{ 1 - \bar r s}{1  - r\bar s}   
                                                             \end{array}\right]    \, ,
$$
from which this convenient formula results:
$$
                                     \frac{p-q}{1  - \bar p q} =  \frac{s-r}{1  - \bar s r} \ ,
$$
provided $p$, $q$, $r$, $s$ are collinear.
\\

\section {Application: relativistic velocities}

\begin{figure}[h]
\centering
\includegraphics[scale=.8]{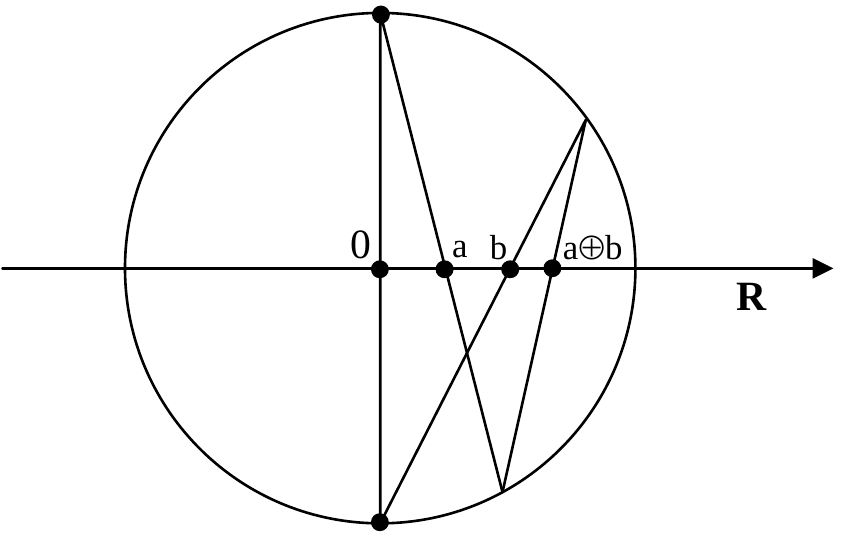}
\caption{\small  Relativistic addition of velocities}
\label{fig:rel}
\end{figure}

Take the real line and the unit circle and as ingredients for the model. 
Consider a quadrilateral that goes through three collinear points, the origin (0), and two points represented 
by real numbers $a$ and $b$.  From Theorem \ref{thm:equiv} it follows that the fourth point on the real line must be
$$
\begin{aligned}
 M(b)M(a)M(0)    
   &=  \left[\begin{array}{rr}   1   & -b \\ 
                                               b  &  -1  \end{array}\right]    
       \left[\begin{array}{rr}   1   & -a \\ 
                                              a    &  -1  \end{array}\right]    
      \left[\begin{array}{rr}    1  &   0 \\ 
                                              0  &  -1  \end{array}\right]    \\[6pt]
&=  \left[\begin{array}{rr}  1+  ab  & -a-b \\ 
                                             a+b  & -1-ab  \end{array}\right]    
 \ \sim \  \left[\begin{array}{cc}   1   & -\frac{a+b}{1+ab} \\ 
                    \frac{a+b}{1+ab}  &  -1  \end{array}\right]    
\end{aligned}
$$
(we used the fact that conjugation does nothing to real numbers).  Thus the fourth point on $\mathbb R$ has 
coordinate 
$$
            a\oplus b   =  \frac{a+b}{1+ab}
$$
But this happens to be the formula for relativistic addition of velocities!  (in the natural units in which the speed of light is 1).
Thus we obtain its geometric interpretation, presented in Figure \ref{fig:rel}.  A conventional derivation of this diagram may be found in \cite{jkd}. 
See also \cite{jkgeo} for an interactive applet.

The segment through the origin does not have to be vertical for the device to work (see Theorem 2.1), 
but is set so for simplicity.

\newpage
\section{Further generalizations}

\begin{figure}[h]
\centering
\includegraphics[scale=.8]{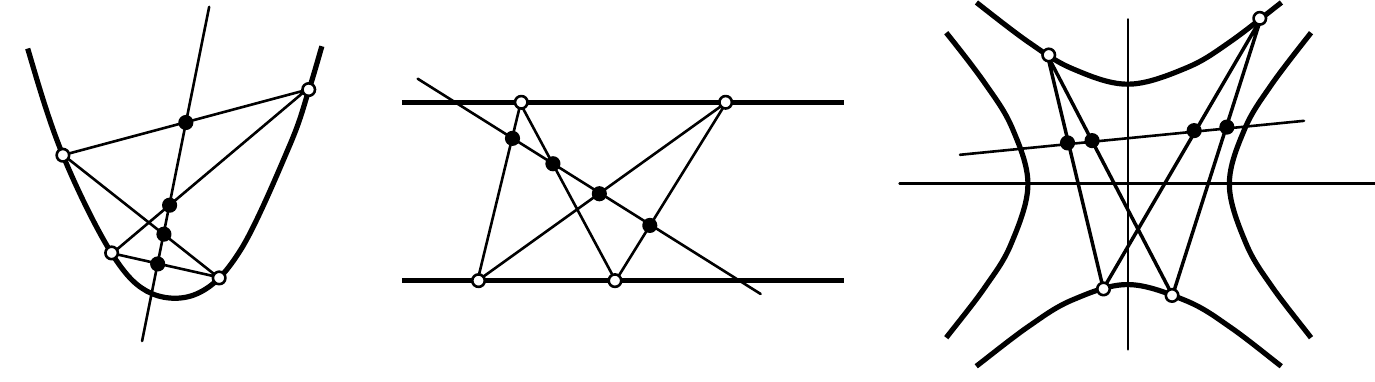}
\caption{\small  Porism on quadrilaterals inscribed in quadrics}
\label{fig:quadrics}
\end{figure}

The porism described here for circles is also valid for any quadric,  
see Figure \ref{fig:quadrics}.

\noindent
What is intriguing, three cases (namely circle, ellipse and a pair of parallel lines) correspond to 
the three possible 2-dimensional Clifford algebras, 
each of the form of ``generalized number plane'' $\{a+be\}$, where $e$ is an ``imaginary unit'' 
whose square is $1$, $-1$ or $0$
(see Table \ref{table}).

\begin{table}[h]
\caption{Algebra-geometry correspondence}
\begin{tabular}{lll}
Algebra                                              &``sphere''             &  formula represented geometrically\\
\hline
$\mathbb{C}$ \ (complex numbers)      &circle                 &    relativistic velocities \\
$\mathbb{D}$ \ (hyperbolic numbers)   &hyperbola         &    trigonometric tangents\\
$\mathbb{G}$ \ (dual numbers)              & parallel lines     &     regular addition
\label{table}
\end{tabular}
\end{table}

The case of the circle corresponds to complex numbers, as described in the previous sections. 
The case of the hyperbola corresponds to ``duplex numbers'', called also hyperbolic numbers or split-complex numbers:
\begin{equation}
\label{eq:D}
                    \mathbb D = \{a+bI \:  : \ a,b\in\mathbb R,\  I^2=1 \}
\end{equation}
They are -- like complex numbers --  a two-dimensional unital algebra except its ``imaginary  unit'' $I$  is 1 when squared.  
As complex numbers are related to rotations, duplex numbers correspond to hyperbolic rotations.  
They were introduced in \cite{Co} as {\it tessarines} and represent a Clifford algebra $\mathbb R_{0,1}$.  
They found useful applications, e.g. in \cite{jkq} a ``hyperbolic quantum mechanics'' was introduced.  
Here is the theorem corresponding to Theorem \ref{thm:moebius}:

\begin{proposition}
\label{prop:D}
The reversion through a point $p\in \mathbb D$  with respect to the hyperbola $x^2 - y^2 = 1$ is 
represented by the matrix
\begin{equation}
\label{eq:moebiusD1}
   \hat p: \ z \quad \longrightarrow  \quad z' = 
                          \left[\begin{array}{rr}   1   & -p \\ 
                                                          \bar p  &  -1  \end{array}\right]   \cdot z 
         =\frac{z-p}{\bar p z - 1}
        \qquad   ( |z|^2=1 )
\end{equation}
while for branch $x^2-y^2 = -1$ it is
\begin{equation}
\label{eq:moebiusD2}
   \hat p: \ z \quad \longrightarrow  \quad z' = 
                          \left[\begin{array}{rr}   -1   & p \\ 
                                                          \bar p  &  1  \end{array}\right]   \cdot z 
         =\frac{-z+p}{\bar p z + 1}
        \qquad   ( |z|^2=-1 )
\end{equation}
\end{proposition}

\noindent
{\bf Proof:}
 For (\ref{eq:moebiusD1})  follow the lines of the proof of the main theorem.  
To get (\ref{eq:moebiusD2}), multiply the ingredients of  (\ref{eq:moebiusD1}) by $I$ to switch the axes, 
use  (\ref{eq:moebiusD1}), and then multiply by $I$ again to restore the original axes. 
\hfill $\square$
\\

For the last case of parallel lines we can use the same matrix calculus as above, 
but replace the algebra by dual numbers
\begin{equation}
\label{eq:G}
                 \mathbb G = \{a+bI\ :\ a,b\in \mathbb R,\, I^2=0 \}
\end{equation}
Transformations  (\ref{eq:moebiusD1}) and (\ref{eq:moebiusD2}) apply for two cases: vertical and horizontal lines, 
respectively ($|z|^2=1$ versus $|z|^2= -1$  in the hyperbolic norm). 
\\

\begin{figure}[h]
\centering
\includegraphics[scale=.8]{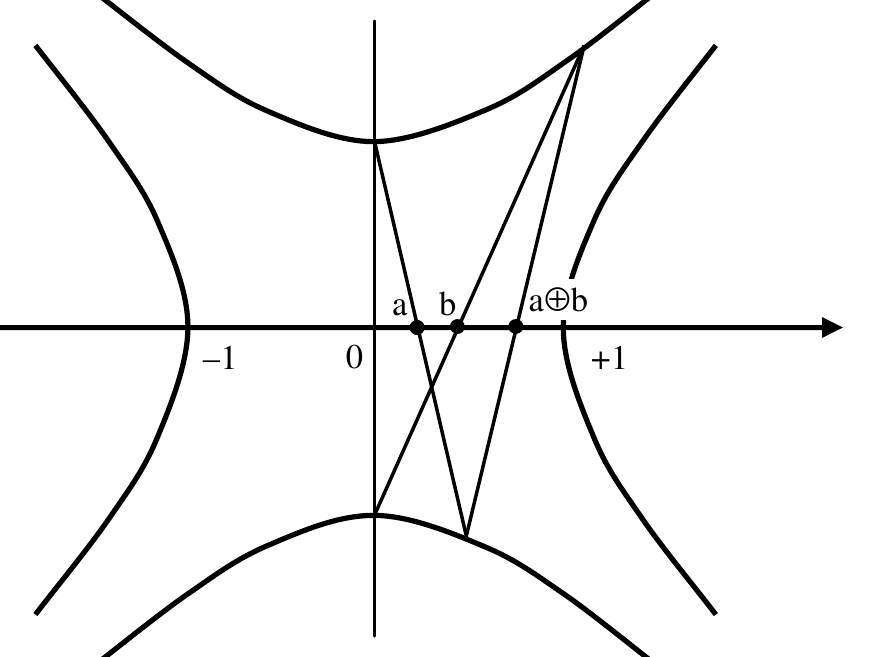}
\caption{\small Trigonometric tangent formula visualized}
\label{fig:trig}
\end{figure}

It is a rather pleasant surprise to find the algebra and geometry interacting at such a basic level in an entirely nontrivial way.  Repeating a construction analogous to the one in Section 5 that results in a geometric tool 
for addition of relativistic velocities to the hyperbolic numbers gives a similar geometric diagram 
representing addition of trigonometric tangents, see Figure \ref{fig:trig}:  
$$
         a\oplus b   = \frac{a+b}{1-ab}
                          = \tan(\arctan  a + \arctan b) 
$$

The final comment concerns the obvious extension of these results implied by inversion of the standard 
configuration like in Figure \ref{fig:porism} through a circle.
M\"obius geometry removes the distinction between circles and lines.
The lines of the quadrilateral and the line $L$  under inversion may become circles.
The point at infinity, where these lines meet, becomes 
under inversion a point that is part of the porism's construction.

\begin{figure}[h]
\centering
\includegraphics[scale=.9]{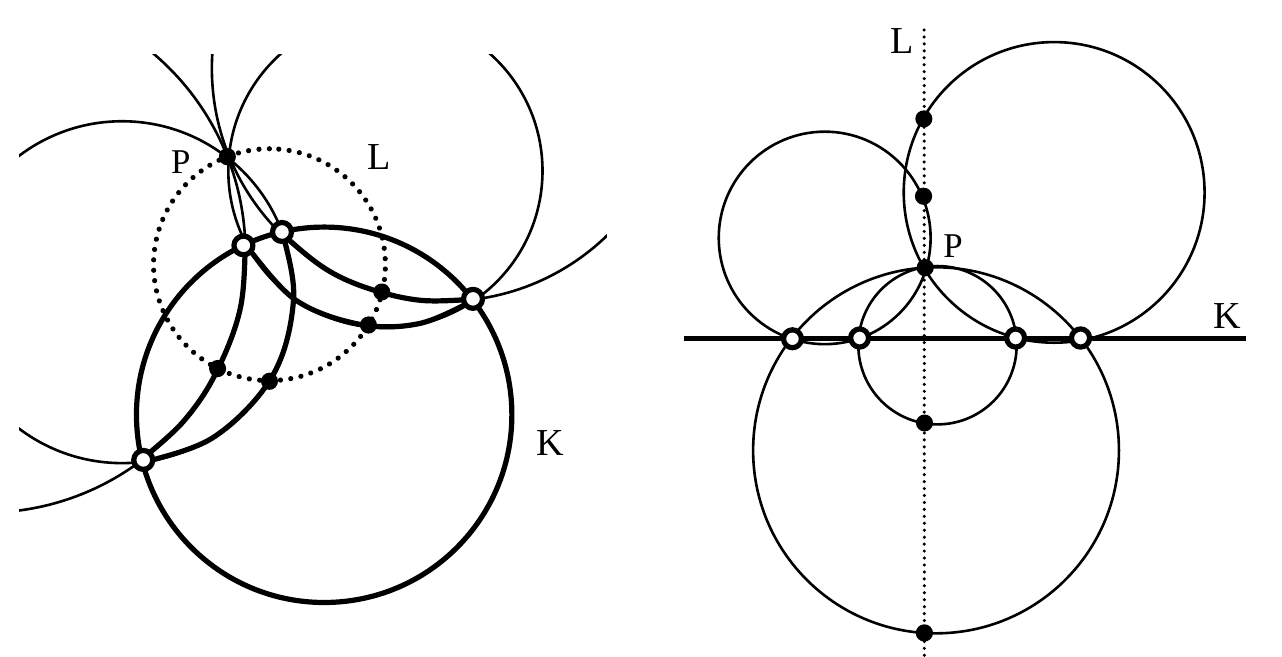}
\caption{\small Porism under inversion extends to circles replacing lines and vice versa}
\label{fig:inverse}
\end{figure}

Thus the general statement is:  given circles $L$ and $K$ and an odd number of 
points $P$, $P_1$, $P_2$, \ldots, $P_{2n}$ on $L$.  Suppose there exist an ordered $2n$-tuple (pencil) of 
circles $C_1,\ldots, C_{2n}$, through $P$ such that $P_i \in C_i$ for every $i$,  and the points of 
intersections $C_i\cap C_{i+1}$ and $C_{2n}\cap C_1$ different from $P$ all lie on the circle $K$, 
then there are infinitely many such 2$n$-tuples of circles. 
Figure \ref{fig:inverse} illustrates the theorem for four circles.  
For interactive version of these and more examples see \cite{jkgeo}.
\\


\section*{Acknowledgments}

The author is grateful to Philip Feinsilver for his interest in this work and his helpful comments.  
Special thanks go also to creators of Cinderella, a wonderful software that allows one to  test quickly
geometric conjectures and to create nice interactive applets.


\end{document}